\documentclass[12pt]{article} 
\usepackage[reqno]{amsmath} 
\usepackage{amssymb, amsthm}
\usepackage{enumerate}  
\usepackage{url} 

\usepackage{hyperref} 
\usepackage{graphicx}

\newtheoremstyle{plainsl}%
	{\topsep}
	{\topsep}
	{\slshape} 
	{}
	{\normalfont\bfseries}
	{.}
	{ }
	{}

\swapnumbers

\theoremstyle{plainsl}
\newtheorem{theorem}{Theorem}[section]
\newtheorem{lemma}[theorem]{Lemma}
\newtheorem{corollary}[theorem]{Corollary}
\newtheorem{proposition}[theorem]{Proposition}

\newcommand\cref[1]{Corollary~\ref{cor:#1}}

\renewcommand\proof{\noindent\textsl{Proof. }}
\newcommand\sqr[2]{{\vbox{\hrule height.#2pt
    \hbox{\vrule width.#2pt height#1pt \kern#1pt
        \vrule width.#2pt}\hrule height.#2pt}}}
\renewcommand\qed{%
	\ifmmode\eqno\sqr53
	\else\nolinebreak\ \hfill\sqr53\medbreak\fi}

\numberwithin{equation}{section}

\newcommand\cx{{\mathbb C}}
\newcommand\compl[1]{{\mkern2mu\overline{\mkern-2mu#1}}}
\newcommand{\Gausswithq}[3]{\left[#1 \atop #2\right]_{#3}}
\newcommand{\ord}{\textrm{ord}}

\title{Perfect state transfer on distance-regular graphs and association schemes}
\author{G.~Coutinho\thanks{Department of Combinatorics and Optimization, University of Waterloo, \texttt{ \{gcoutinho, cgodsil\}@uwaterloo.ca}}\ \thanks{Supported by Capes Foundation, Ministry of Education, Brazil.} , C.~Godsil\footnotemark[1] \footnotemark[3], K.~Guo\thanks{Supported in part by NSERC.} \thanks{ Department of Mathematics, Simon Fraser University, \texttt{krystalg@sfu.ca}} , and F.~Vanhove \thanks{Supported by the Research Foundation Flanders-Belgium (FWO-Vlaanderen).}\ \ \thanks{Department of Mathematics, Ghent University.}}

\begin{document}
\maketitle

\begin{center}
\em In memory of Fr\'{e}d\'{e}ric Vanhove, died November 27, 2013. \end{center}

\begin{abstract}
We consider the representation of a continuous-time quantum walk in a graph $X$ by the matrix $\exp(-it A(X))$. We provide necessary and sufficient criteria for distance-regular graphs and, more generally, for graphs in association schemes to have perfect state transfer. Using these conditions, we provide several new examples of perfect state transfer in simple graphs.
\end{abstract}

\section{Introduction} 

For  given graph $X$ with adjacency matrix $A$, the transition matrix of the quantum walk  at time $t$ is given by $U = e^{i t A}$. If $u$ and $v$ are distinct vertices of graph $X$, we say that $X$ admits \textsl{perfect state transfer} if there is a time $\tau$ such that
\[U(\tau) e_u = \lambda e_v\]
for some $\lambda \in \mathbb C$. 

The main problem we are concerned about is that of determining which graphs admit perfect state transfer. This problem is solved for paths and hypercubes (see \cite{PhysRevLett.92.187902}) and for circulant graphs (see \cite{Petkovi2011300}). It was also extensively discussed for cubelike graphs in general (see \cite{MR2504947} and \cite{MR2811131}). The effect of certain graph operations in perfect state transfer has also been considered in \cite{MR2655511} and in \cite{MR2933057}. A generalization to signed graphs is considered in \cite{MR3076338}. Some recent surveys are \cite{Kendon:2011-03-01T00:00:00:1546-1955:422} and \cite{MR2852516}.

In general, perfect state transfer is understood to be a relatively rare phenomenon. One of the first infinite families of graphs admitting perfect state transfer was given in \cite{godsil11}. In the known examples of graphs admitting perfect state transfer, the times at which perfect state transfer occur tend to be at $\frac{\pi}{2}$ or $\frac{\pi}{4}$. The known classes of graphs where perfect state transfer occurs at a time earlier than $\frac{\pi}{4}$ are the complete bipartite graphs $K_{2,n}$, examples from \cite{chan2013} and various examples given in this paper, including the Hadamard graphs. 

In this paper, we will examine perfect state transfer in distance-regular graphs and, more generally, in graphs contained in association schemes. This problem was already studied in \cite{PhysRevA.77.022315} and in \cite{MR2992400}, but here we will present a necessary and sufficient condition that can be easily tested.  We determine precisely which graphs in the known infinite families of distance regular graphs admit perfect state transfer.  In addition, for all graphs listed in \cite[Chapter 14]{MR1002568}, we determine which admit perfect state transfer. 

\section{Quantum walks} 

Quantum walk is an important concept in quantum algorithms. A quantum walk is a quantum analog of a classical random walk. Quantum walk algorithms have been studied and shown to perform exponentially or polynomially better for various black box problems. We refer to \cite{Childs09, Ambainis07, MR1638221, CCDFGS03} for further background on such results.

There are several ways of defining a quantum analog to a random walk. We focus on the \textsl{continuous-time quantum walk}, generated by the adjacency matrix of a graph. Let $X$ be a graph with adjacency matrix $A$. The continuous-time quantum walk at time $t$ is given by the unitary operator  $e^{-itA}$ and takes place in the Hilbert space whose elementary basis corresponds to the vertices of $X$. In particular, we address the problem of when $X$ admits a transfer of state between two vertices without a loss of information with respect to a continuous-time quantum walk. We call this phenomenon \textsl{perfect state transfer} and we formalize it as follows.

For every non-negative real number $t$, we recall that
\[U_A(t) = e^{i t A} = \sum_{k \geq 0} \frac{(i t)^k}{k!} A^k.\]
We omit the subscript $A$ whenever the context is clear. Observe that
\[\overline{U(t)} = U(-t).\]
and operator $\overline{U(t)}$ is the continuous-time quantum walk on $X$.  The matrix $A$ is symmetric, hence $U(t) = U(t)^T$. Therefore
\[U(t)^*   U(t) = U(-t)  U(t) = I,\]
so $U(t)$ is a unitary operator. For notational convenience, we will work with $U(t)$, but all the results we present are also true for $U(-t)$.

If $u$ and $v$ are distinct vertices of a graph $X$, recall that $X$ admits \textsl{perfect state transfer} if there is a time $\tau$ such that
\[U(\tau) e_u = \lambda e_v\]
for some $\lambda \in \mathbb C$. Note that $|\lambda| = 1$ because $U(t)$ is unitary for all $t$. Note also that perfect state transfer from $u$ to $v$ implies it from $v$ to $u$. For this reason, we typically refer to $uv$-pst in a graph $X$. This definition was first considered by Christandl \textsl{et al.} in \cite{PhysRevLett.92.187902}.

We also say that $X$ is \textsl{periodic at $u$} at time $\tau \neq 0$ if
\[U(\tau) e_u = \lambda e_u\]
for some $\lambda \in \mathbb C$. Note that if $X$ admits $uv$-pst at time $\tau$, then $X$ is periodic at both $u$ and $v$ at time $2 \tau$. A graph is called \textsl{periodic} if it is periodic at every vertex. It was proved by Kay \cite[Section D]{Kay2011} that if $uv$-pst takes place in a graph $X$, neither of the vertices can be involved in perfect state transfer with a third vertex.

Suppose that the adjacency matrix $A$ of a graph has eigenvalues $\theta_0 > \theta_1 > ... > \theta_n$. The matrix $A$ is real and symmetric, so it is orthogonally diagonalizable and therefore admits a spectral decomposition:
\[A = \sum_{r=0}^n \theta_r E_r.\]
Since $U_A(t)$ is defined in terms of powers of $A$, the following holds:
\[U(t) = \sum_{r=0}^n e^{it \theta_r} E_r.\]
For this reason, the spectrum of a graph can give valuable information about the existence of perfect state transfer.

\section{Association schemes and distance-regular graphs}

A set of $n\times n$ matrices with $0$-$1$ entries $\{A_0,...,A_d\}$ is a (symmetric) \textsl{association scheme} if the following properties hold:
\begin{enumerate}[(1)]
\item the identity matrix $I$ belongs to the set, say $A_0 = I$;
\item $A_i$ is symmetric for $i = 0, \ldots, d$;
\item $\displaystyle \sum_{k=0}^d A_k = J$; and,
\item there exist integers $p_{ij}^k$, said to be \textsl{intersection numbers of the scheme}, such that:
\[A_i A_j = \sum_{k=0}^d p_{ij}^k A_k.\]
\end{enumerate} 

The matrices $\{A_0,...,A_d\}$ are said to be \textsl{classes} of the scheme. Note that (2) and (4) together imply that these matrices pairwise commute. The matrix algebra spanned by these matrices is the \textsl{Bose-Mesner algebra} of the association scheme, usually denoted by the symbol $\mathcal A$. It is a commutative algebra, and the fact that $J$ belongs to the algebra implies that any matrix in the algebra has a constant row and column sum. The classical reference for the theory of association schemes is \cite{MR882540}.

The Schur product of matrices $M$ and $N$ is defined as:
\[(M \circ N)_{ab} = M_{ab} \cdot N_{ab}.\]
Note that:
\[A_i \circ A_j = \left\{ \begin{array}{ll}
A_i & \text{ if $i = j$,} \\ 0 & \text{ otherwise.}
\end{array}  \right.\]
The matrices $\{A_0,...,A_d\}$ are idempotents with respect to the Schur product, and so the Bose-Mesner algebra is closed under Schur product. A set of symmetric and pairwise commuting matrices can be simultaneously diagonalized, hence there exist idempotents for the conventional matrix product $\{E_0,...,E_m\}$ which also form a basis for the algebra. In particular $m = d$ and any matrix in $\mathcal A$ has at most $d+1$ different eigenvalues.

Using both bases for the Bose-Mesner algebra of an association scheme, we can define coefficients $P_{ji}$ and $Q_{ij}$ as follows: 
\begin{align}
A_i & = \sum_{j=0}^d P_{ji} E_j , \label{Pij} \\
E_j & = \frac{1}{n} \sum_{i=0}^d Q_{ij} A_i .\label{Qij} 
\end{align}
The matrix $P$ with entries $P_{ji}$ is called the \textsl{eigenmatrix} of the scheme.

The construction of non-trivial association schemes is not an easy task. We now describe a construction based on graphs that exhibit a high level of regularity.

A connected graph is said to be \textsl{distance-regular} if there exist numbers $b_i,\ c_i$ for $i \geq 0$ such that for any two vertices $u$ and $v$ at distance $i$, the number of neighbours of $v$ at distance $i-1$ from $u$ is $c_i$ and the number of neighbours of $v$ at distance $i+1$ from $u$ is $b_i$ (note that these numbers do not depend on the choice of $u$ and $v$). This definition implies the existence of numbers $a_i$ which are the number of neighbours of $v$ at distance $i$ from $u$. They do not depend on $u$ and $v$ because $b_0$ is the valency of the graph and $a_i$ is determined in terms of $b_i$ and $c_i$ by the following equation:
\begin{align} \label{eq1} b_0 = c_i + a_i + b_i. \end{align}
A distance-regular graph  is said to have \textit{classical parameters} $(d,b,\alpha,\beta)$ where $\alpha,\beta \in \mathbb{R}$ and $b \in \mathbb{R} \setminus \{ 0 \}$ if
\begin{eqnarray}
b_i&=&\biggl(\Gausswithq{d}{1}{b}-\Gausswithq{i}{1}{b}\biggr)\biggl(\beta-\alpha\Gausswithq{i}{1}{b}\biggr), \nonumber \\
c_i&=&\Gausswithq{i}{1}{b}\biggl(1+\alpha\Gausswithq{i-1}{1}{b}\biggr), \nonumber
\end{eqnarray}
for every $i \in \{ 0,1,\ldots,d \}$. Here, $\Gausswithq{i}{1}{b}:=i$ if $b=1$ and $\Gausswithq{i}{1}{b}:= \frac{b^i-1}{b-1}$ if $b \neq 1$.

We say that a graph \textsl{belongs} to an association scheme if its adjacency matrix is contained in the Bose-Mesner algebra of the scheme. The results in this paper are motivated by a search for a characterization of perfect state transfer in distance-regular graphs, but they generalize to graphs belonging to association schemes. Unlike the examples cited at the end of the introduction, the class of distance-regular graphs is a large class of graphs of which a full characterization is unlikely to ever be made. The classical reference is \cite{MR1002568}. We follow \cite{MR1220704}  and introduce preliminary concepts and definitions.

Suppose $X$ is a distance-regular graph of diameter $d$. The list of parameters $\{b_0,b_1,...,b_{d-1} ;\ c_1,c_2,...,c_d\}$ is called the \textsl{intersection array} of the graph. The following proposition provides a well-known necessary condition for an array of numbers to be the intersection array of a distance-regular graph.

\begin{proposition}\cite[Proposition 4.1.6]{MR1002568} \label{intersecarray}
If $\{b_0,b_1,...,b_{d-1} ;\ c_1,c_2,...,c_d\}$ is the intersection array of a distance-regular graph, then $b_0 > b_1 \geq b_2 \geq ... \geq b_{d-1} > 0$ and $1 = c_1 \leq c_2 \leq ... \leq c_d \leq b_0$.
\end{proposition}

We define the graphs $X_i$ as the graphs with vertex set $V(X)$ and two vertices adjacent if and only if they are at distance $i$ in $X$. Let $A = A(X)$ and define distance matrices $A_i(X) = A(X_i)$. Define $A_0$ as the identity matrix. By definition of a distance-regular graph, the matrices $\{A_0,\ A,\ A_2,\ ...,\ A_d\}$ satisfy:
\begin{align*}AA_i = b_{i-1} A_{i-1} + a_i A_i + c_{i+1}A_{i+1}. \tag{1}\end{align*}
In particular, the matrices $A_i$ can be written as a polynomial of degree $i$ in $A$ for all $i \geq 1$. 

By induction and using the equation above, one can prove that there exist \textsl{intersection numbers of the graph} $p_{ij}^k$, which depend only on the intersection array, such that:
\[A_i A_j = \sum_{k = 0}^d p_{ij}^k A_k.\]

Then, if $X$ is a distance-regular graph of diameter $d$, the matrices $\{A_0,A_1,\\ ...,A_d\}$ form an association scheme, and the intersection numbers of the graph coincide with those of the scheme.

An association scheme with intersection numbers $p_{ij}^k$ is said to be \textsl{metric} if, for a particular ordering of the relations and for all $i$, $j$ and $k$, we have:
\[p_{ij}^k \neq 0\  \Rightarrow \  i+j \geq k\]
and
\[p_{ij}^{i+j} \neq 0\]
when $i+j\leq d$. An association scheme with a given ordering of relations is \textsl{$P$-polynomial} if the entries of the $k$-th row of the eigenmatrix $P$ are evaluations of a polynomial $p_k(x)$ of degree $k$ at the eigenvalues of $A_1$.

Not all association schemes arise from distance matrices of distance-regular graphs. Those which do are characterized by the following classical theorem. 

\begin{theorem} \cite[Theorems 5.6 and 5.16]{MR0384310}
Given an association scheme $\{A_0,...,A_d\}$, the following conditions are equivalent:
\begin{enumerate}[(i)]
\item The scheme (with the given ordering) is metric.
\item The scheme is $P$-polynomial.
\item There exists a distance-regular graph $X$ of diameter $d$ such that $A(X_i) = A_i$.
\end{enumerate}
\end{theorem}

A distance-regular graph $X$ of diameter $d$ is said to be \textsl{primitive} if the graphs $X_i, i\in\{1,\ldots,d\}$ are all connected, and \textsl{imprimitive} otherwise. If the valency of  an imprimitive distance-regular graph $X$ is at least $3$, then $X$ is either bipartite or $X_d$ is the disjoint union of cliques of the same size, case in which the graph is said to be \textsl{antipodal} and the cliques in $X_d$ are said to be \textsl{fibres} of the graph (see for instance \cite[Theorem 4.2.1]{MR1002568}).

\section{Characterization}

First, we present a necessary condition for a graph belonging to an association scheme to admit perfect state transfer, which originally appeared in \cite[Theorem 4.1 and Lemma 6.1]{godsil11}. Here, we provide an explicit proof.

\begin{theorem}\cite[Theorem 4.1 and Lemma 6.1]{godsil11} \label{pst implis ant size 2}
Let $X$ be a graph that belongs to an association scheme with $d$ classes and with adjacency matrix $A = A(X)$. If $X$ admits perfect state transfer at time $\tau$, then there is a permutation matrix $T$ with no fixed points and of order two such that $U_A(\tau) = \lambda T$ for some $\lambda \in \cx$. Moreover, $T$ is a class of the scheme. If the graph is distance-regular of diameter $d$, then it must be antipodal with fibres of size $2$ and $T = A_d$.
\end{theorem}

\proof
From 
\[U(t) = \sum_{r=0}^n e^{it \theta_r} E_r,\]
we see that $U(t)$ belongs to the Bose-Mesner algebra $\mathcal A$ of the scheme. If, for some $\tau$,
\[U(\tau)_{u,v} = \lambda \quad \text{ with } \quad |\lambda| = 1,\]
then consider the $0$-$1$ matrix $T$ which is the unique element of the basis of Schur idempotents with $1$ in the $(u,v)$ position. We have that
\[
U(\tau) = \sum_{i=0}^d \alpha_i A_i
\]
with $\alpha_i \in \cx$.  So the coefficient of $T$ is $\lambda$. Since $U(\tau)$ has exactly one non-zero entry in the $u$th row and matrices in the scheme have constant row sum, the coefficients of $A_i \neq T$ are all equal to 0. Then $U(\tau) = \lambda T$. 

Observe that the $u$th row of $T$ is the basis vector $e_v$. The row  and column sum of $T$ are constant and must then be equal to $1$. Therefore $T$ is a permutation matrix. Since $T \neq A_0$, the permutation represented by $T$ has no fixed points. Since $U(\tau)$ is symmetric, we have that $T$ is a permutation matrix of order 2. 

Now if $X$ is distance-regular, let $i$ be the index for which $A_i = T$. Suppose $0 < i < d$, and so $d>1$. Hence $v$ is the only vertex at distance $i$ from $u$, then $b_{i-1} = 1$. By Proposition \ref{intersecarray}, this implies that $b_j = 1$ for all $j \geq i$. In particular, there will be a unique vertex at distance $d$ from $u$, and this vertex will have degree $1$, hence the graph is $K_2$ and $d=1$, a contradiction. Therefore $A_d = T$, and the distance-regular graph is antipodal with fibres of size $2$.
\qed

We may say more about the spectrum of graphs in association schemes admitting perfect state transfer. Recall that a graph $X$ is periodic if there exists a time $t$ such that $U_A(t)$ is a diagonal matrix.

\begin{lemma}\label{lem:nec} Let $X$ be a graph belonging to an association scheme. If perfect state transfer occurs on $X$, then all its eigenvalues are integers.
 \end{lemma}
\proof
By Theorem \ref{pst implis ant size 2}, $X$ will be periodic at time $2 \tau$. Since $X$ is regular, $X$ has at least one integral eigenvalue. Then, by Corollary 3.3 of \cite{godsil11}, all eigenvalues of $X$ must be integers.
\qed

Consider an association scheme with $d$ classes and minimal idempotents $\{E_0,...,E_d\}$ that contains a permutation matrix $T$ of order two as a relation. Note that the eigenvalues of $T$ are either $+1$ or $-1$, and so $TE_j = \pm E_j$. We define a partition $(\mathcal{I}_T^+,\mathcal{I}_T^-)$ of $\{0,...,d\}$ as $j \in \mathcal{I}_T^+$ if $TE_j = E_j$, and $j \in \mathcal{I}_T^-$ otherwise. We might drop the subscript $T$ if the context is clear.

If $x \in \mathbb{Z}$, denote by $\ord_2(x)$ the exponent of $2$ in the factorization of $x$ (here we convention $\ord_2(0) = + \infty$). We now present our main theorem.

\begin{theorem}\label{thm:drgpst} Let $X$ be a graph with distinct eigenvalues $\lambda_0 > ... > \lambda_d$ belonging to an association scheme. Let
\[ \alpha := \gcd(\{\lambda_0-\lambda_k : k \in \{0,...,d\}\}).\]
Perfect state transfer occurs in $X$ at time $t$ between vertices $u$ and $v$ if and only if all of the following hold:
\begin{enumerate}[(i)]
\item The relation of the scheme that contains $(u,v)$, say $T$, is a permutation matrix of order 2 with no fixed points;
\item the eigenvalues of $X$ are integers; and,
\item $\ord_2(\lambda_0-\lambda_j) > \ord_2(\alpha)$ for all $j \in \mathcal{I}^+$ and $\ord_2(\lambda_0-\lambda_\ell) = \ord_2(\alpha)$ for all $\ell \in \mathcal{I}^-$. Moreover, $t$ is an odd multiple of $\frac{\pi}{\alpha}$.
\end{enumerate} 
\end{theorem}

\proof By Theorem \ref{pst implis ant size 2} and Lemma \ref{lem:nec}, conditions (i) and (ii) are necessary for $X$ to have perfect state transfer at time $t$. Suppose conditions (i) and (ii) hold.  

We may write 
\[ U(t) = \sum_{j\in \mathcal{I}^+} e^{i \lambda_j t } E_{j} +  \sum_{\ell \in \mathcal{I}^-} e^{i \lambda_{\ell} t} E_{\ell}
\]
and, for any $\gamma \in \mathbb{C}$,
\[ \gamma T = \sum_{j \in \mathcal{I}^+} \gamma E_{j}  + \sum_{\ell \in \mathcal{I}^-} - \gamma E_{\ell}.\]
Since $U(t)$ and $T$ are both elements of the Bose-Mesner algebra of the scheme and the idempotents $E_j$ form a basis for this algebra, we see that $U(t) = \gamma T$ for some $t$ and some $\gamma$ if and only if the corresponding coefficients of every $E_{j}$ and $E_{\ell}$ are equal. It follows from the Perron-Frobenius theory (see for instance \cite[Theorem  2.2.1]{MR2882891}) that $0 \in \mathcal{I}^+$, so $\gamma = e^{i \lambda_0 t}$. Thus, for all $j \in \mathcal{I}^+$, $e^{i \lambda_j t} = \gamma$ if and only if $t(\lambda_0-\lambda_j)$ is a multiple of $2 \pi$; and, for all $\ell \in \mathcal{I}^-$, $e^{i \lambda_\ell t} = -\gamma$ if and only if $t(\lambda_0-\lambda_{\ell})$ is an odd multiple of $\pi$. A time $t$ satisfying these conditions can be chosen if and only if condition (iii) holds. Moreover, the minimum time $t$ satisfying these conditions is $\tau = \tfrac{\pi}{\alpha}$. Perfect state transfer happens at another time $\tau_1$ if and only if $\tau_1$ is an odd multiple of $\tau$.
\qed

Given the $P$-matrix of an association scheme, it is easy to see which eigenvalues of $X$ have idempotents corresponding to the eigenvalues $1$ or $-1$ of $T$. Note that there may be more than one perfect matching among the classes of the scheme, in which case condition (iii) must be tested in each one of them. In this case, at most one of them can satisfy condition (iii). In the case of distance-regular graphs, we may find the partition $(\mathcal{I}^+, \mathcal{I}^-)$ explicitly using the following lemma.

\begin{lemma}\cite[Proposition 11.6.2]{MR2882891}\label{lem:lambdamu} Let $X$ be a distance-regular graph of diameter $d$. Suppose $X$ is antipodal with classes of size two and let $\lambda_0 > \ldots > \lambda_d$ be the distinct eigenvalues of $X$ with corresponding idempotents $E_0,\ldots , E_d$. Then
\[
A_d E_j = (-1)^j E_j.
\] 
\end{lemma}

Using Lemma \ref{lem:lambdamu}, we now simplify our main theorem in the case where $X$ is distance-regular.

\begin{corollary} \label{cor:drgpst} Let $X$ be a distance-regular graph of diameter $d$ and distinct eigenvalues $\lambda_0 > \ldots > \lambda_d$. As before, let
\[ \alpha := \gcd(\{\lambda_0-\lambda_k : k \in \{0,...,d\}\}).\]
Perfect state transfer occurs on $X$ at time $t$ if and only if all of the following hold
\begin{enumerate}[(i)]
\item $X$ is antipodal with classes of size 2;
\item the eigenvalues of $X$ are integers;
\item \label{numbertheoretic} $\ord_2(\lambda_0-\lambda_j) > \ord_2(\alpha)$ if $j$ is even and $\ord_2(\lambda_0-\lambda_\ell) = \ord_2(\alpha)$ if $\ell$ is odd. Moreover, $t$ is an odd multiple of $\frac{\pi}{\alpha}$.
\end{enumerate}
\end{corollary} 
\proof
For condition (i), see Theorem \ref{pst implis ant size 2}. Condition (ii) is the same as in Theorem \ref{thm:drgpst}, whereas condition (iii) of this corollary is condition (iii) of Theorem \ref{thm:drgpst} rewritten using Lemma \ref{lem:lambdamu}.
\qed

\section{Perfect state transfer in distance-regular graphs}

In this section we examine examples of perfect state transfer in distance-regular graphs. The distance-regular graphs of diameter $2$ are also known as \textsl{(connected) strongly regular graphs} (see \cite[section 1.3]{MR1002568} or \cite[Chapter 10]{MR1829620} as a reference for what follows below).

The parameters of a strongly regular graphs are usually presented as $(v,k,a,c)$, where $v$ is the number of vertices, $k$ is the valency, $a$ is the number of common neighbours of two adjacent vertices, and $c$ is the number of common neighbours of two non-adjacent vertices. Such a graph is antipodal if and only if being at distance 0 or 2 is an equivalence relation, and in terms of the parameters, this implies that $c=k$ and $a=2k-v$. It follows that the graphs are complete multipartite with classes of size $v-k$. By Corollary \ref{cor:drgpst}, perfect state transfer happens only if $v-k =2$, case in which the graph is the complement of a disjoint union of $m$ copies of $K_2$. The distinct eigenvalues of such graphs are
\[\{2m - 2,\ 0,\ -2\}.\]
Hence, using Corollary \ref{cor:drgpst}, we have:
\begin{corollary}
Perfect state transfer happens in distance-regular graphs of diameter $2$ if and only if the graph is the complement of a disjoint union of $m$ copies of $K_2$ with $m$ even. In that case, perfect state transfer happens at time $\frac{\pi}{2}$.
\end{corollary}

\subsection{Distance-regular graphs of diameter 3}

We now study in detail distance-regular graphs of diameter 3. Perfect state transfer can occur only if the graphs are antipodal with classes of size $2$, case in which the graphs are also called \textsl{Taylor graphs}, and are equivalent to combinatorial objects known as \textsl{regular two-graphs}. We will completely determine when perfect state transfer occurs among Taylor graphs, and by doing that, we will provide a large list of examples of simple graphs admitting perfect state transfer, all of these examples so far unknown to the best of our knowledge.

Suppose $X$ is a graph and $\pi$ is a partition of the vertex set of $X$ satisfying the following two conditions:
\begin{enumerate}[(1)]
\item each cell is an independent set, and
\item between any two cells either there are no edges or there is a perfect matching.
\end{enumerate}
Let $Y = X / \pi$ be the simple graph that has the cells of $\pi$ as its vertices, and two vertices adjacent if and only if there is a matching in $X$ between the corresponding cells. Then we say that $X$ is a \textsl{covering} of $Y$. If $Y$ is connected, condition (2) implies that the cells of $\pi$ have the same size, say $r$, and in this case we say that $X$ is an \textsl{$r$-fold covering} of $Y$.

Antipodal distance-regular graphs of diameter $3$ must be covering graphs of the complete graph (see \cite[section 4.2B]{MR1002568} or see \cite[Theorem 2.1]{Godsil1992205}). The parameters of the intersection array of a distance-regular $r$-fold covering of $K_n$ will depend only on $n$, $r$ and $c$, where $c$ stands for the number of common neighbours shared by vertices at distance $2$ of the covering (see \cite[Lemma 3.1]{Godsil1992205}). The intersection array of such graphs is:
\[\{n-1,\ (r-1)c,\ 1 ;\ 1,\ c,\ n-1 \},\]
and then we say the graph is a distance-regular $(n,r,c)$ covering of $K_n$. Based on the intersection array, we can compute the eigenvalues in terms of $n$, $r$ and $c$ (see, for instance, \cite[section 3]{Godsil1992205}). Define the parameters $\delta = n-rc-2$ and $\Delta = \delta^2 + 4(n-1)$. The set of distinct eigenvalues is:
\begin{align}\{ n-1,\ \tfrac{\delta + \sqrt{\Delta}}{2}  ,\ -1,\ \tfrac{\delta - \sqrt{\Delta}}{2}\}. \label{eq:spectrum}\end{align}

The following theorem will permit us to organize our search.

\begin{theorem}\cite[Theorem 3.6]{Godsil1992205} \label{thm:4.2}
For fixed $r$ and $\delta$, there are only finitely many feasible parameter sets for distance-regular covers of $K_n$, unless $\delta = -2$, $\delta = 0$ or $\delta = 2$.
\end{theorem}

Note that $r$ is the size of the antipodal classes. By Corollary \ref{cor:drgpst}, we only have  to consider the case when $r =2$. Distance-regular $2$-fold coverings of the complete graph are equivalent to strongly regular graphs with parameters $\left(n-1,\ n-c-2,\ n-\tfrac{3c}{2} -3,\ \tfrac{n-c-2}{2}\right)$. These are the graphs induced in the neighbourhood of every vertex of the covering (see \cite[Chapter 11]{MR1829620}).

We now begin our classification of graphs in this class where perfect state transfer occurs. Let $X$ be a distance-regular $(n,2,c)$ covering of the complete graph. We proceed via the case analysis suggested by Theorem \ref{thm:4.2}.

\begin{theorem} \label{nothappensdelta0}
A distance-regular cover of $K_n$ with $\delta = 0$ does not admit perfect state transfer.
\end{theorem}
\proof
First note that we are already assuming $r=2$. If $\delta = 0$, then $n = 2c + 2$. In order for perfect state transfer to occur using the conditions of Corollary \ref{cor:drgpst}, the eigenvalues must be integers, hence $\Delta = 4(n-1)$ must be a square. So $n-1$ must be a square, hence $n$, which is even, is congruent to $2$ modulo $4$. Note that $n = (n-1) - (-1)$ and $(n-1)$ and $(-1)$ are eigenvalues with the same parity. Therefore condition (\ref{numbertheoretic}) of Corollary \ref{cor:drgpst} implies that if perfect state transfer occurs in this case, it will occur at time $\frac{\pi}{b}$, $b$ an odd number, and so the differences between eigenvalues with different parities must be odd. But $\sqrt{n-1}$ and $n-1$ are both odd, hence their difference is even. We conclude that in this case perfect state transfer cannot occur.
\qed

For $\delta \neq 0$, the following proposition says that the other cases occur in pairs (see for instance \cite[p.431]{MR1002568}).
\begin{proposition}\label{pairsofcovers}
Suppose $X$ is a distance-regular $2$-fold cover of $K_n$. Then the distance 2 graph $X_2$ is a distance-regular $2$-fold cover of $K_n$, with $\delta(X_2) = -\delta(X)$.
\end{proposition}

We now determine the structure of the covers for which $\delta = -2$. We will give a simplified form of \cite[Lemma 8.2 ]{Godsil1992205}. A \textsl{Hadamard matrix} of order $n$ is a $n\times n$ matrix $H$ with entries in $\{1,-1\}$ such that 
\[ HH^T = H^TH = I_{n}.
\]
Let $X$ be a $(n,2,c)$ distance-regular cover with $\delta = -2$. Order the antipodal pairs from $1$ to $n$ and for each antipodal pair, define an arbitrary ordering of its vertices. Define a square matrix $B$ of order $n$ indexed by the antipodal pairs of $X$ as follows:
\[
B_{ij} := \begin{cases} 0 & \text{if } i = j; \\
+1 & \text{if the matching between pairs } $i$ \text{ and } $j$ \text{ agrees with their } \\
& \text{respective orderings}; \\
-1 & \text{otherwise.} \end{cases} \]
\begin{theorem} \label{thm:hadcover}
The matrix $(B+I)$ is a symmetric Hadamard matrix with constant diagonal. Conversely, every symmetric Hadamard matrix with constant diagonal and order $n$ yields a $\left( n,2,\dfrac{n}{2} \right)$ distance-regular cover. 
\end{theorem}

If $X$ is a $(n,2,c)$ cover with $\delta = -2$, it follows from arithmetic conditions on the multiplicities of the eigenvalues of $X$ that $n$ must be a square (see \cite[Lemma 3.7]{Godsil1992205}). This can also be derived as a necessary condition for the existence of symmetric Hadamard matrices of constant diagonal.

It is also a known fact that Hadamard matrices can only exist when $n$ is 1, 2 or a multiple of $4$ (see, for instance, \cite[section 1.8]{MR1002568}).

\begin{theorem} \label{hadamard3}
Every $2$-fold distance-regular cover of $K_n$ with $\delta = -2$ admits perfect state transfer at time $\frac{\pi}{\sqrt{n}}$. For $\delta = +2$, perfect state transfer occurs if and only if $n$ is divisible by $8$, and in that case it occurs at time $\frac{\pi}{2}$.
\end{theorem}
\proof
Note that $n$ must be an even square of an integer.

Let $X$ be a $2$-fold distance-regular cover of $K_n$ with $\delta = -2$. Using \ref{eq:spectrum} in this case, its set of distinct eigenvalues will be:
\[ \{n-1,\ \sqrt{n}-1,\ -1,\ -\sqrt{n}-1\}. \]
Applying Corollary \ref{cor:drgpst}, we see that perfect state transfer will occur at time $\dfrac{\pi}{\sqrt{n}}$. If $\delta = +2$, then the set of distinct eigenvalues is:
\[ \{n-1,\ \sqrt{n}+1,\ -1,\ -\sqrt{n}+1\}. \]
If $ord_2(n) = 2$, then $n$ and $n - \sqrt{n}-2$ are both divisible by $4$. If $ord_2(n)>2$, then perfect state transfer will occur at time $\tfrac{\pi}{2}$.
\qed

Constructions for symmetric Hadamard matrices with constant diagonal are known for all $n$ which are a power of $4$ (see \cite{MR2455527}). Perfect state transfer was already known for the case $\delta = -2$ (see \cite{godsil11}), but unknown for the case $\delta = +2$.

Now we move to the case where $\delta \notin \{0,-2,2\}$. As we saw in Proposition \ref{pairsofcovers}, for every cover with parameter $\delta$, there exists a corresponding cover with parameter $-\delta$.
\begin{theorem} \label{pstdistd3}
Let $X$ be a $(n,2,c)$ distance-regular cover, and let $X_2$ be the corresponding $(n,2,n-2-c)$ cover. Suppose $\delta = n-2c-2 \notin \{0,-2,2\}$. If $\delta \equiv 2 \bmod 4$, then perfect state transfer occurs either in $X$ or in $X_2$ at time $\frac{\pi }{2a}$, for some $a$ which is an odd integer. If $\delta$ is odd or a multiple of $4$, perfect state transfer does not occur in either $X$ or $X_2$.
\end{theorem}
\proof
Let $\rho>0$ and $\sigma<0$ be the eigenvalues of $X$ which are neither $n-1$ nor $-1$. Let $m_\rho$ and $m_\sigma$ be their corresponding multiplicities, which can be computed in terms of $n$, $r$ and $c$ (see \cite[section 3]{Godsil1992205}). Then:
\[m_\rho - m_\sigma = \frac{n \delta}{\sqrt{\Delta}},\]
where $\Delta = \delta^2 + 4(n-1)$. This difference must be an integer, and $\delta \neq 0$ implies that $\Delta$ must be a perfect square. Note that $\sigma = \frac{1}{2}(\delta - \sqrt{\Delta})$ is an algebraic integer, hence an integer. Thus $ \sqrt{\Delta} - \delta$ is even. Suppose that $\Delta = (2t + \delta)^2$. The parameter $n$ is now parametrized as $n = 1+t^2 + t \delta$.

Note that if $n$ is odd, then perfect state transfer cannot occur by Corollary \ref{cor:drgpst}. If $\delta$ is odd, then $n = \delta + 2c + 2$ is odd. If $t$ is even, then $n$ is also odd. So suppose $\delta$ is even and $t$ is odd.

If $\delta$ is a multiple of $4$, then $n = 1 + t^2 + t\delta \equiv 2 \bmod 4$. Perfect state transfer occurs only if $(n-1) - \rho$ and $(n-1) - \sigma$ are both odd. Hence $\rho$ and $\sigma$ are both even. But note that $\sigma = -t$ and $t$ is odd.

If $\delta \equiv 2 \bmod 4$, then $n = 1 + t^2 + t\delta \equiv 0 \bmod 4$. Let $\rho_2>0$ and $\sigma_2>0$ be the non-trivial eigenvalues of $X_2$. Note that:
\[ \rho = \delta+t \quad \text{and} \quad \sigma = -t,\]
and
\[ \rho_2 = t \quad \text{and} \quad \sigma_2 = -\delta-t.\]

If $t \equiv 3 \bmod 4$, then $(n-1) - (\delta + t)$ and $(n-1) +t$ are both congruent to $2$ modulo $4$, so perfect state transfer occurs in $X$ at time $\frac{\pi}{2a}$, where $2a$ is the greatest common divisor of the differences of the eigenvalues. If $t \equiv 1 \bmod 4$, then $(n-1) - t$ and $(n-1) - (-\delta - t)$ are both congruent to $2$ modulo $4$, so perfect state transfer occurs in $X_2$ at time $\frac{\pi}{2a}$, where $2a$ is the greatest common divisor of the differences of the eigenvalues.
\qed

Note that in the case where $\delta \equiv 2 \bmod 4$, the theorem above implies that perfect state transfer occurs at least in one of $X$ and $X_2$ and it is possible that perfect state transfer occurs in both.

Table 1 below contains parameters sets for which a corresponding cover admits perfect state transfer, for $n < 280$. The strongly regular graph on the first neighbourhood of a vertex in each of these covers is given in the rightmost column. We consult Andries Brouwer's website (http://www.win.tue.nl/ $\sim$aeb/graphs/srg/srgtab.html) to either provide a construction for such a graph, hence a construction for the cover, or to state that no such construction is known.

\begin{table}[!h]\caption{pst in distance-regular graphs of diameter 3} \label{table} \ \\
$\begin{array}{c|c|c|c|c} 
n & c & \delta & \text{time} & \text{construction of strongly regular graph} \\ \hline
28&10&6&{\pi}/{2}& \text{complement of Schl\"afli graph} \\ 
76&42&-10&{\pi}/{2}& \text{not known} \\ 
96&40&14&{\pi}/{4}& \text{not known} \\ 
96&54&-14&{\pi}/{6}& \text{not known} \\ 
120&54&10&{\pi}/{6}& \text{complement of $O^-(8,2)$ polar graph} \\ 
136&70&-6&{\pi}/{2}& \text{complement of $O^+(8,2)$ polar graph} \\ 
148&66&14&{\pi}/{2}& \text{not known} \\ 
176&72&30&{\pi}/{4}& \text{complement of the one below} \\ 
176&102&-30&{\pi}/{2}& \text{line graph of Hoffman-Singleton graph} \\ 
244&130&-18&{\pi}/{2}& \text{not known} \\ 
276&162&-50&{\pi}/{6}& \text{complement of McLaughlin graph McL.2 / U4(3).2} \\ 
\end{array}$
\end{table}

\subsection{Diameter larger than 3}

We present examples of distance-regular graphs with diameter larger than $3$ admitting perfect state transfer. In this subsection we do not provide a full characterization as we did for the diameter $3$ case, but we discard certain infinite families of distance-regular graphs with classical parameters.

A \textsl{Hadamard graph} of order $n\geq 2$ is a graph $X(H)$ obtained from a $n \times n$ Hadamard matrix $H$ as follows: $X(H)$ has a pair of vertices $\{r^+, r^-\}$ for each row and a pair of vertices $\{c^+, c^-\}$ for each column. If the entry of $H$ in row $r$ and column $c$ is $1$, then $X(H)$ has edges $\{c^-, r^-\}$ and $\{c^+, r^+\}$. If the entry of $H$ in row $r$ and column $c$ is $-1$, then $X(H)$ has edges $\{c^-, r^+\}$ and $\{c^+, r^-\}$. These graphs are not to be confused with the covers of Theorem \ref{thm:hadcover}. The Hadamard graphs of order $n$ are $n$-regular graphs with $4n$ vertices. They are distance-regular of diameter 4, antipodal with classes of size two, bipartite and have intersection array $\{n, n-1, \frac{n}{2}, 1; 1, \frac{n}{2}, n-1, n\}$. Their distinct eigenvalues are
\[\{n,\ \sqrt{n},\ 0,\ -\sqrt{n},\ -n\},\]
and so, by Corollary \ref{cor:drgpst}, perfect state transfer happens if and only if $n$ is a perfect square, and at time $\frac{\pi}{\sqrt{n}}$.
 
We organize the example above and the remaining examples of this section in the following corollary. We refer to \cite[Chapters 9, 11 and 13]{MR1002568} for information about these graphs.

\begin{corollary} \label{distd4more}
The graphs below are examples of distance-regular graphs of diameter larger than three admitting perfect state transfer.
\begin{enumerate}[(i)]
\item Hamming $d$-cubes. Number of vertices: $2^d$. Valency: $d$. Diameter: $d$. Distinct eigenvalues: $\{d-2i : i = 0,...,d\}$. Time of perfect state transfer: $\frac{\pi}{2}$.
\item Halved $2d$-cubes. Number of vertices: $2^{2d-1}$. Valency: $\binom{2d}{2}$. Diameter: $d$. Distinct eigenvalues: $\left\{\binom{2d}{2} - 2i(2d-i) : i=0,...,d\right\}$. Time of perfect state transfer: $\frac{\pi}{2}$.
\item Hadamard graphs of order $n$ if and only if $n$ is a perfect square. Exists for infinitely many values of $n$, in particular for all $n$ which are powers of $4$. Number of vertices: $4n$. Valency: $n$. Diameter: $4$. Time of perfect state transfer: $\frac{\pi}{\sqrt{n}}$.
\item Meixner graph  (\cite[Example 3.5]{MMW07}). Number of vertices: $1344$. Valency: $176$. Diameter: $4$. Distinct eigenvalues: $\{176,44,8,-4,-16\}$. Time of perfect state transfer: $\frac{\pi}{12}$.
\item Coset graph of the once shortened and once truncated binary Golay code. Number of vertices: $1024$. Valency: $21$. Diameter: $6$. Distinct eigenvalues: $\{21,9,5,1,-3,-7,-11\}$. Time of perfect state transfer: $\frac{\pi}{4}$.
\item Coset graph of the shortened binary Golay code. Number of vertices: $2048$. Valency: $22$. Diameter: $6$. Distinct eigenvalues: $\{22, 8, 6, 0, -2,\\ -8, -10\}$. Time of perfect state transfer: $\frac{\pi}{2}$.
\item Double coset graph of truncated binary Golay code. Number of vertices: $2048$. Valency: $22$. Diameter: $7$. Distinct eigenvalues: $\{22, 10, 6, 2, -2,\\ -6, -10, -22\}$. Time of perfect state transfer: $\frac{\pi}{4}$.
\item Double coset graph of binary Golay code. Number of vertices: $4096$. Valency: $23$. Diameter: $7$. Distinct eigenvalues: $\{23,9,7,1,-1,-7,-9,\\-23\}$. Time of perfect state transfer: $\frac{\pi}{2}$.
\end{enumerate} 
\end{corollary}

\begin{corollary}
No graph in the following infinite families of antipodal distance-regular graphs with classes of size two admits perfect state transfer.
\begin{enumerate}[(i)]
\item Johnson graphs $J(2n,n)$ for $n>1$. Their distinct eigenvalues are $\{(n-j)^2-j : j \in \{0,...,n\} \}$.
\item Doubled Odd graphs on $2n+1$ points. Their distinct eigenvalues are $\{(-1)^j(n+1-j) : j \in \{0,1,...,n-1,n,n+2,n+3,...,2n+2\}\}$.
\end{enumerate}
\end{corollary}

We also checked graphs with diameter larger than 3 depicted in \cite[Chapter 14]{MR1002568} that do not belong to the infinite families above.

\begin{corollary}\label{cor:notpst}
None of the following antipodal distance-regular graphs with classes of size two and diameter larger than three admit perfect state transfer.
\begin{enumerate}[(i)]
\item Wells graph of diameter 4. Some eigenvalues are not integral.
\item Double Hoffman-Singleton Graph of diameter 5. Distinct eigenvalues are $\{7,3,2,-2,-3,-7\}$.
\item Double Gewirtz Graph of diameter 5. Distinct eigenvalues are $\{10,4,2,\\-2,-4, -10\}$.
\item Double 77-Graph of diameter 5. Distinct eigenvalues are $\{16,6,2,-2,\\-6,-16\}$.
\item Double Higman-Sims Graph of diameter 5. Distinct eigenvalues are $\{22,8,2,-2,-8,-22\}$.
\item Dodecahedron of diameter 5. Distinct eigenvalues are not integer.
\end{enumerate}
\end{corollary}

\section{Perfect state transfer in graphs belonging to association schemes}

In this section, we present examples of perfect state transfer in graphs which are not distance-regular but which belong to association schemes. In the first subsection, we show how the tensor product can be used to construct association schemes that admit perfect state transfer. In the second subsection, we exhibit examples of perfect state transfer in sporadic association schemes. In \cite{chan2013}, Chan gives examples of graphs admitting perfect state transfer which are constructed from taking unions of distance classes in the binary Hamming scheme. 

\subsection{Bipartite doubles}

Consider two association schemes $\{A_0,...,A_d\}$ and $\{B_0,...,B_e\}$. Consider the set of matrices obtained by taking the tensor product of the matrices in both schemes:
\[\{A_i \otimes B_j : \ 0 \leq i \leq d \text{ and } 0 \leq j \leq e \}. \]
This is an association scheme with $(de + d + e)$ classes (see \cite[Chapter 3]{MR2047311}). 

If $A$ and $B$ are the adjacency matrices of graphs $X$ and $Y$, then the matrix $A \otimes B$ is the adjacency matrix of the graph on $V(X) \times V(Y)$ where $(x_1,y_1)$ is adjacent to $(x_2,y_2)$ if and only if $x_1$ is adjacent to $x_2$ in $X$ and $y_1$ to $y_2$ in $Y$. This is said to be the \textsl{direct product} of $X$ and $Y$ (note that the graph obtained does not depend on the order of the product). In the literature, the direct product is also referred to as the weak product or the tensor product of graphs. 

In the context of this paper, we present examples of direct products of a distance-regular graph $X$ with $K_2$ that admit perfect state transfer. Such construction is also referred to as the \textsl{bipartite double} of $X$. Since these graphs belong to association schemes, the conditions of Theorem \ref{thm:drgpst} are sufficient to determine if perfect state transfer occurs. Other examples of perfect state transfer in direct products of graphs were studied in \cite{MR2934678}. Some of these graphs were already considered in Corollary \ref{cor:notpst}. 

\begin{theorem}
Let $V(K_2) = \{v_1,v_2\}$. Suppose $X$ is distance-regular on $n$ vertices with eigenvalues $\theta_0 > ... > \theta_d$, and let $\theta_i = 2^{f_i} m_i$, where $m_i$ is an odd integer. For any vertex $u \in V(X)$, $X \times K_2$ admits perfect state transfer between $(u,v_1)$ and $(u,v_2)$ if and only if both conditions hold:
\begin{enumerate}[(i)]
\item For all $i$, $f_i = a$ for some constant $a$.
\item For all $i$ and $j$, $m_i \equiv m_j \bmod 4$.
\end{enumerate}
Under these conditions, perfect state transfer occurs at time $\frac{\pi}{2\gcd\{\theta_0,...,\theta_d\}}$.
\end{theorem}
\proof
It follows from the definition of direct product that the eigenvalues of $X \times K_2$ are $\pm \theta_i$. Note that the matrix $I_n \otimes A(K_2)$ is a permutation matrix of order two and with no fixed points belonging to the scheme, and that
\[(I_n \otimes A(K_2) )(E_{+ \theta_i}) = E_{+ \theta_i} \quad \text{and} \quad (I_n \otimes A(K_2)) (E_{- \theta_i}) = - E_{- \theta_i}.\]
Condition (iii) of Theorem \ref{thm:drgpst} uses $\lambda_0$ as a pivot, but in this context it is equivalent to the two conditions below. Here $\alpha$ is the gcd of the differences of the eigenvalues of $X \times K_2$.
\begin{enumerate}[(1)]
\item $\ord_2(\theta_i - \theta_j) > \ord_2 (\alpha)$ for all $i$ and $j$; and
\item $\ord_2(\theta_i - (-\theta_j)) = \ord_2 (\alpha)$ for all $i$ and $j$.
\end{enumerate}
If $i = j$, then (2) is equivalent to Condition (i) of the statement. Under Condition (i), (1) and (2) reduce to 
\begin{enumerate}
\item[(1')] $\frac{m_i - m_j}{2}$ even, for all $i$ and $j$; and
\item[(2')] $\frac{m_i + m_j}{2}$ odd, for all $i$ and $j$.
\end{enumerate}
which is equivalent to Condition (ii). The time follows from the expression for time given in Theorem \ref{thm:drgpst}.
\qed

The following corollaries exhibit some new examples of perfect state transfer. Note that the distance between vertices involved in perfect state transfer in these examples is the odd girth of the graph. The definition of these graphs can be found in \cite[Chapter 10]{MR1829620} for strongly regular graphs or in \cite{MR1002568} for distance-regular graphs. For generalized quadrangles, a more detailed account can be found in \cite{payne2009finite}.

\begin{corollary} \label{bipdoublesrg}
The following bipartite doubles of strongly regular graphs admit perfect state transfer. 
\begin{enumerate}[(i)]

\item Bipartite doubles of the point graphs of generalized quadrangles $GQ(s,t)$ whenever $q$ is a prime power and one of the following conditions hold:
\begin{itemize}
\item[] $s=q-1,t=q+1$ with $q\equiv 0 \bmod4$,
\item[] $s=q,t=q^2$ with $q\equiv 7\bmod 8$,
\item[] $s=q^3,t=q^2$ with $q\equiv 7 \bmod 8$,
\item[] $t=1$ , with $s\equiv 7 \bmod 8$.
\end{itemize}
Number of vertices: $2(st+1)(s+1)$. Valency: $s(t+1)$. Perfect state transfer at $\frac{\pi}{4}$.

\item Bipartite doubles of the complements of the point graphs of generalized quadrangles $GQ(s,t)$ whenever $q$ is a prime power and one of the following conditions hold:
\begin{itemize}
\item[] $s=q,t=q^2,$ with $q\equiv 3 \bmod 4$: perfect state transfer at $\frac{\pi}{2 q}$,
\item[] $s=q^2,t=q$, with $q \equiv 3 \bmod 4$: perfect state transfer at $\frac{\pi}{2 q}$,
\item[] $s=q-1,t=q+1$, with $q$ even : perfect state transfer at $\frac{\pi}{2 }$,
\item[] $s=q+1,t=q-1$, with $q$ even: perfect state transfer at $\frac{\pi}{2 }$,
\item[] $s=q^2,t=q^3$, with $q\equiv 3 \bmod 4$: perfect state transfer at $\frac{\pi}{2 q^2}$,
\item[] $s=q^3,t=q^2$, with $q\equiv 3 \bmod 4$: perfect state transfer at $\frac{\pi}{2 q^2}$.
\end{itemize}
Number of vertices: $2(st+1)(s+1)$. Valency: $s^2t$.

\item Bipartite doubles of orthogonal array graphs $OA(n,m)$ if $\ord_2(n)\geq \ord_2(m)+2$. Number of vertices: $2n^2$. Valency: $m(n-1)$. Perfect state transfer at $\frac{2 \pi}{gcd(n,4 m)}$.

\item Bipartite doubles of complements of orthogonal array graphs $OA(n,m)$ if $\ord_2(n)\geq \ord_2(m-1)+2$. Number of vertices: $2n^2$. Valency: $n^2-m(n-1)-1$. Perfect state transfer at $\frac{2 \pi}{gcd(n,4 (m-1))}$.
\end{enumerate}
\end{corollary}

We now note that the eigenvalues of a distance-regular graph with classical parameters $(d,b,\alpha,\beta)$ are given by the following formula (see \cite[Corollary 8.4.2]{MR1002568}):
\[\Gausswithq{d-j}{1}{b}\left(\beta-\alpha\Gausswithq{j}{1}{b}\right)-\Gausswithq{j}{1}{b},j\in\{0,\ldots,d\}.\]

\begin{corollary} \label{bipdoubledrg}
The following bipartite doubles of distance-regular graphs with classical parameters admit perfect state transfer.
\begin{enumerate}[(i)]
\item Bipartite doubles of Grassmann graphs $J_q(n,d)$ ($n\geq 2 d$) for $n$ even, $d$ odd and $q \equiv 3 \bmod 4$. Classical parameters: $\left(d,q,q,q\frac{q^{n-d}-1}{q-1}\right)$. Perfect state transfer at $\frac{\pi}{2}$.
\item Bipartite doubles of Hamming graphs $H(d,q)$ when $\ord_2(q) \geq \ord_2(4d)$. Classical parameters: $(d,1,0,q-1)$. Perfect state transfer at $\frac{2\pi}{\gcd(q,4d)}$.
\item Bipartite doubles of Doob graphs of odd diameter. Classical parameters: $(d,1,0,3)$. Perfect state transfer at $\frac{\pi}{2}$.
\item Bipartite doubles of unitary dual polar graphs $^2A_{2d-1}(q)$ and $^2A_{2d}(q)$, both cases when $\ord_2(q+1) \geq \ord_2 (4d)$. Classical parameters: $(d,q^2,0,q)$ and $(d,q^2,0,q^3)$, respectively. Perfect state transfer at $\frac{2\pi}{\gcd(q+1,4d)}$.
\item Bipartite doubles of parabolic and symplectic dual polar graphs $B_d(q)$ and $C_d(q)$ when $q \equiv 3 \bmod 4$ and $d$ is odd. Classical parameters: $(d,q,0, q)$. Perfect state transfer at $\frac{\pi}{2}$.
\item Bipartite doubles of half dual polar graph of diameter $d$ on $2d$-spaces when $\ord_2(q+1) \geq \ord_2(4d)$. Classical parameters: $\left(d,q^2,q^2+q,q\frac{q^{2 d-1}-1}{q-1}\right)$. Perfect state transfer at $\frac{\pi}{2 gcd(d,q+1)}$.
\item Bipartite doubles of half dual polar graph of diameter $d$ on $2 d+1$-spaces  when $\ord_2(q+1)(q^2+1))\geq \ord_2(4 d)$. Classical parameters: $\left(d,q^2,q^2+q ,q\frac{q^{2d+1}-1}{q-1}\right)$. Perfect state transfer at
\[ \frac{2 \pi}{gcd\left((q+1)(q^2+1),4\frac{q^{2 d+1}-1}{q^2-1}\frac{q^{2  d}-1}{q-1}\right)}.\]
\item Bipartite doubles of exceptional graphs of Lie type when $q\equiv  11 \bmod 12$ or when $q\equiv 3,7 \bmod 12$. Classical parameters: $\left(3,q^4,q\frac{q^4-1}{q-1},q\frac{q^9-1}{q-1}\right)$. In the first case, perfect state transfer at time $\frac{\pi}{6}$. In the second, at time $\frac{\pi}{2}$.
\item Bipartite doubles of all affine $E_6$ graphs when $q$ is even. Classical parameters: $(3,q^4,q^4-1,q^9-1)$. Perfect state transfer at $\frac{\pi}{2}$.
\item Bipartite doubles of all alternating forms graphs when $q$ is even. Classical parameters: $(d,q^2,q^2-1,q^{2 d-1}-1)$ and $(d,q^2,q^2-1,q^{2 d+1}-1)$ for forms on $2d$- and on $2d+1$-spaces, respectively. Perfect state transfer at $\frac{\pi}{2}$.
\item Bipartite doubles of all Hermitian forms graphs when $q$ is even. Classical parameters: $(d,-q,-q-1,-(-q)^d-1)$. Perfect state transfer at $\frac{\pi}{2}$.
\end{enumerate}
\end{corollary}

\subsection{Sporadic association schemes}

Here we present examples of perfect state transfer coming from sporadic association schemes. A quick description of these examples can be found after the corollary.

\begin{corollary}\label{pstsporadic} 
The following association schemes contain at least one graph admitting perfect state transfer.
\begin{enumerate}[(i)]
\item the second relation of the $Q$-bipartite $4$-class scheme based on $GQ(s,s^2)$ described in \cite{MR2739499} has perfect state transfer at time $\frac{\pi}{2 s}$ for all odd $s$;
\item the first relation in the extended $Q$-bipartite double of the $3$-class scheme arising from a linked system of $l$ symmetric $(16 s^2, 2 s ( 4 s -1), 2 s (2 s -1))$ designs has perfect state transfer at time $\frac{\pi}{4 s}$ for all odd $l$;
\item the orthogonality relation on the tight spherical $7$-design consisting of the $240$ vectors of the $E_8$ root system has perfect state transfer at time $\frac{\pi}{6}$;
\item the relation corresponding to inner product $\frac{1}{3}$ on the tight spherical $7$-design consisting of a specific set of $4600$ vectors derived from the Leech lattice has perfect state transfer at time $\frac{\pi}{6}$;
\item in general, in any scheme obtained from a tight spherical $7$-designs in $\mathbb{R}^d$ with $d=3 w^2-4$ for some integer $w\geq 2$, the first relation has perfect state transfer at time $\frac{\pi}{w(w-1)}$ if $w\equiv 3\bmod{4}$, the second at $\frac{\pi}{2(w^2-1)}$  if $w$ is even, and the third at $\frac{\pi}{w(w+1)}$ if $w$ is even;
\item the orthogonality relation on the unique (up to isomorphism) tight spherical $11$-design in $\mathbb{R}^{24}$ has perfect state transfer at time $\frac{\pi}{18}$.
\end{enumerate}
\end{corollary}

\subsubsection*{Penttila and Williford's $4$-class scheme}

Penttila and Williford \cite{MR2739499} considered a $Q$-bipartite $4$-class scheme, based on a $GQ(s,s^2)$ with a subquadrangle $GQ(s,s)$ for $s>2$.  Its $P$-matrix is given by: 
\[P= \left( \begin {array}{ccccc} 1& \left( s-1 \right)  \left( {s}^{2}+1
 \right) & \left( {s}^{2}-2\,s \right)  \left( {s}^{2}+1 \right) &
 \left( s-1 \right)  \left( {s}^{2}+1 \right) &1\\ \noalign{\medskip}1
&{s}^{2}+1&0&-{s}^{2}-1&-1\\ \noalign{\medskip}1&s-1&-2\,s&s-1&1
\\ \noalign{\medskip}1&-s+1&0&s-1&-1\\ \noalign{\medskip}1&- \left( s-
1 \right) ^{2}&2\,s \left( s-2 \right) &- \left( s-1 \right) ^{2}&1
\end {array} \right).\]
The second non-identity relation corresponds to points in the $GQ(s,s^2)$, not in the subquadrangle, being non-collinear, different from each other's antipode and not collinear to each other's antipode either.

\subsubsection*{Linked systems of symmetric designs} 
Martin, Muzychuk and Williford \cite[Theorem 3.1]{MMW07} proved a method to construct $Q$-bipartite association schemes with a matching, the ``extended $Q$-bipartite double", from certain $Q$-polynomial association schemes.  In \cite[Theorem 3.6]{MMW07} they apply this to a linked system of $l+1$ symmetric  $(16 s^2, 2 s ( 4 s -1), 2 s (2 s -1))$ designs with $l\geq 2$.  The $P$-matrix of the derived $Q$-bipartite $4$-class scheme is then given by:
\[ P=\left( \begin {array}{ccccc} 1&16\,l{s}^{2}&32\,{s}^{2}-2&16\,l{s}^{2
}&1\\ \noalign{\medskip}1&4\,ls&0&-4\,ls&-1\\ \noalign{\medskip}1&0&-2
&0&1\\ \noalign{\medskip}1&-4\,s&0&4\,s&-1\\ \noalign{\medskip}1&-16\,
{s}^{2}&32\,{s}^{2}-2&-16\,{s}^{2}&1\end {array} \right).\]

\subsubsection*{Schemes from tight spherical designs}
Consider a finite set $X$ of unit vectors in Euclidean space $\mathbb{R}^d, d\geq 1$.  Its \textsl{degree} $s$ is the number of values assumed by the inner product between distinct vectors in $X$.  We say $X$ is \textsl{antipodal} if $\forall x\in X: -x \in X$.  Delsarte,Goethals and Seidel \cite[Theorem 6.8]{MR0485471} proved that if $X$ is antipodal, then $|X|\leq 2 \binom{d+s-2}{d-1}$.   In case of equality, we call $X$ a tight $(2s-1)$-design.  The $s$ inner products between distinct vectors in $X$ then define a $Q$-polynomial $s$-class association scheme on $X$ by \cite[Theorem 7.5]{MR519045}, the parameters of which can be computed using \cite[Theorem 3.6 and Remark 7.6]{MR519045}, and in particular only depend on $d$ and $s$.  Bannai and Damerell \cite{MR519045,MR576179} proved that tight $(2s-1)$-designs can only exist for $s\in\{1,2,3,4,6\}$ and that $s=6$ implies that $d=24$.

If $s=4$, then $d=3 w^2-4$ for some integer $w\geq 2$ (see for instance \cite[pp.1401-1402]{MR2535394}).  The $P$-matrix is given by:
\[P= \left( \begin {array}{ccccc} 1&\frac{{w}^{4} \left( 3\,{w}^{2}-5
 \right) }{2} & \left(\! -\!16{w}^{2}+6{w}^{4}+10 \right)\! \left(\!-\!1+{w}
^{2} \right) & \frac{{w}^{4} \left( 3{w}^{2}-5 \right)}{2} &1
\\ \noalign{\medskip}1& \frac{ \left(3{w}^{2}-5\right) {w}^{3}}{2}&0&-
 \frac{\left( 3{w}^{2}-5 \right) {w}^{3}}{2}&-1\\ \noalign{\medskip}1&
 \left( {w}^{2}-2 \right) {w}^{2}&-2\, \left( -1+{w}^{2} \right) ^{
2}& \left( {w}^{2}-2 \right) {w}^{2}&1\\ \noalign{\medskip}1&-w&0&
w&-1\\ \noalign{\medskip}1&-{w}^{2}&-2+2\,{w}^{2}&-{w}^{2}&1
\end {array} \right),\]
where the relations correspond to inner products $1,\frac{1}{w},0,-\frac{1}{w}$ and $-1$, respectively.  Bannai and Sloane \cite{MR617634} proved that if $d=8$, the only tight $7$-design consists of the $240$ vectors of the $E_8$ root system (where the orthogonality relation has perfect state transfer at time $\frac{\pi}{6}$), and that if $d=23$, the only tight $7$-design is a specific set of $4600$ vectors derived from the Leech lattice (where the relation corresponding to inner product $\frac{1}{3}$ has perfect state transfer at time $\frac{\pi}{6}$).
 
If $s=6$, then by Bannai and Sloane \cite{MR617634} there is only one tight $11$-design in $\mathbb{R}^{24}$, consisting of the $196560$ minimal vectors of the Leech lattice.  The $P$-matrix is given by:
\[P= \left( \begin {array}{ccccccc} 1&4600&47104&93150&47104&4600&1
\\ \noalign{\medskip}1&2300&11776&0&-11776&-2300&-1
\\ \noalign{\medskip}1&1000&1024&-4050&1024&1000&1
\\ \noalign{\medskip}1&350&-704&0&704&-350&-1\\ \noalign{\medskip}1&76
&-320&486&-320&76&1\\ \noalign{\medskip}1&-10&16&0&-16&10&-1
\\ \noalign{\medskip}1&-20&64&-90&64&-20&1\end {array} \right),\]
where the relations correspond to inner products $1,\frac{1}{2},\frac{1}{4},0,-\frac{1}{4},-\frac{1}{2},-1$, respectively.

We refer to the survey by Bannai and Bannai \cite{MR2535394} for further information on spherical designs.

\section{Another approach using the eigenvalues of the mixing matrix.}

Given a graph $X$ with adjacency matrix $A = A(X)$, we define the \textsl{mixing matrix} $M_A(t) = M(t)$ as follows: 
\[
M(t) = U(t) \circ \overline{U(t)}
\]
It is a symmetric, doubly-stochastic matrix with non-negative real entries. If $X$ is distance-regular and admits perfect state transfer at time $\tau$, then, by \ref{pst implis ant size 2}, $X$ is antipodal with antipodal classes of size 2 and the mixing matrix is a permutation matrix. In fact, $M(t) = A_d$. On the other hand, if $M(\tau)$ is a permutation matrix of order two and no fixed points, then clearly $X$ admits perfect state transfer at time $\tau$.

In general, the matrix $M(t)$ always belongs to the Bose-Mesner algebra of the association scheme determined by the distance matrices of $X$. Then, we can find an explicit formula for the eigenvalues of $M(t)$ in terms of the $P$ and $Q$ matrices of the association scheme, and use this to investigate the existence of perfect state transfer in $X$.

\begin{lemma} If $X$ is a distance-regular graph with $n$ vertices, diameter $d$ and $M(t)$ is its mixing matrix, then 
\begin{equation}\label{eigM}
M(t) E_{\ell} = \left( \sum_{k = 0}^d \frac{1}{n^2} \left(\sum_{r = 0}^d e^{i \theta_r t} Q_{kr}  \right)\left(\sum_{j = 0}^d \overline{e^{i \theta_j t}} Q_{kj}  \right) P_{\ell k} \right) E_{\ell}
\end{equation}
for $\ell = 0, \ldots, d$. 
\end{lemma}

\proof Recall the definition of the $P$ and $Q$ matrices of a distance regular graph $X$. Then, we may write $U(t)$ using the basis formed by the classes of the scheme:
\[ \begin{split}
U(t) &= \sum_{r=0}^d e^{i\theta_r t} E_r \\
&=  \sum_{r=0}^d \frac{e^{i\theta_r t}}{n} \left( \sum_{k=0}^d Q_{kr} A_k \right) \\
&= \sum_{k=0}^d  \left( \sum_{r=0}^d \frac{e^{i\theta_r t}Q_{kr}}{n} \right) A_k .
\end{split}
\]
Then we may compute $U(t) \circ \compl{U(t)}$ as follows:
\[
U(t) \circ \compl{U(t)} = \sum_{k=0}^d \left( \sum_{r=0}^d \frac{e^{i\theta_r t}Q_{kr}}{n}  \right) \left(\sum_{j=0}^d \frac{\compl{e^{i\theta_j t}}Q_{kj}}{n}  \right) A_k.
\]
We may apply another change of basis to obtain the eigenvalues of $U(t) \circ \compl{U(t)}$ as follows:
\[ \begin{split}
U(t) \circ \compl{U(t)} &= \sum_{k=0}^d \frac{1}{n^2} \left( \sum_{r=0}^d e^{i\theta_r t}Q_{kr}  \right) \left(\sum_{j=0}^d \compl{e^{i\theta_j t}}Q_{kj} \right) \left( \sum_{\ell = 0}^d P_{\ell k} E_\ell  \right) \\
&=  \sum_{\ell = 0}^d \left(\sum_{k=0}^d \frac{1}{n^2} \left( \sum_{r=0}^d e^{i\theta_r t}Q_{kr}  \right) \left(\sum_{j=0}^d \compl{e^{i\theta_j t}}Q_{kj} \right) P_{\ell k}  \right) E_{\ell}.
\end{split}
\]
Then the eigenvalues of $U(t) \circ \compl{U(t)}$  are 
\[  \sum_{k=0}^d \frac{1}{n^2} \left( \sum_{r=0}^d e^{i\theta_r t}Q_{kr}  \right) \left(\sum_{j=0}^d \compl{e^{i\theta_j t}}Q_{kj} \right) P_{\ell k}
\]
as claimed. \qed 

The following lemma allow us to use the eigenvalues computed above to determine if the distance-regular graph has perfect state transfer between pairs of antipodal vertices.

\begin{lemma} If $M$ is a $n\times n$ matrix such that 
\begin{itemize}
\item $M = M^T$;
\item $M$ has row and column sums equal to 1;
\item $M$ has non-negative entries; and,
\item the eigenvalues of $M$ are $1$ and $-1$,
\end{itemize}
then $M$ is a permutation matrix. \end{lemma}

\proof Consider the $i$th row of $M$. The row sum is 
\[\sum_{k = 1}^n M_{ik} = 1. \]
Since $MM^T = M^2$, and $M^2 = I$, we have that 
\[\sum_{k = 1}^n (M_{ik})^2 = 1. \]
Because $0 \leq M_{ik} \leq 1$ for all $k$, we have
\[1 = \sum_{k=1}^n (M_{ik})^2 \leq \sum_{k=1}^n M_{ik} = 1.\]
So equality holds in the inequality, and that happens if and only if $(M_{ik})^2 =M_{ik}$ for all $k$. That is equivalent to $M_{ij} = 1$ for some $j \in \{1,...,n\}$, and $0$ otherwise.
\qed 

Now we have the following theorem using the previous lemmas along with Lemma \ref{lem:lambdamu}. 

\begin{theorem}
Consider a distance-regular graph $X$ with mixing matrix $M(t)$.  Perfect state transfer occurs at time $t$ if and only if 
\[
(-1)^{\ell} = \sum_{k = 0}^d \frac{1}{n^2} \left(\sum_{r = 0}^d e^{i \theta_r t} Q_{kr}  \right)\left(\sum_{j = 0}^d \overline{e^{i \theta_j t}} Q_{kj}  \right) P_{\ell k}
\]
for $\ell = 0, \ldots, d$. 
\end{theorem} 

\proof Since $M(t) = U(t) \circ \overline{U(t)}$, we have that $M(t)$ is symmetric, has row and column sums equal to $1$ and is non-negative. Then, by the previous lemma, if $M(t)$ has all eigenvalues in $\{1, -1\}$, then $M(t)$ is a permutation matrix and perfect state transfer occurs. Then, $M(t) = A_d$ and so must have the same eigenvalues with respect to the idempotents $E_{\ell}$ for $\ell = 0, \ldots, d$, which are given by Lemma \ref{lem:lambdamu}. Conversely, if $X$ has perfect state transfer from $u$ to $v$ at time $t$, then $M(t) = A_d$ is a permutation matrix of order 2 and 
\[ M(t) E_{\ell} = A_d E_{\ell}.
\] \qed

\section{Open problems}

Cubelike graphs are examples of graphs belonging to an association scheme. Although work has been done in \cite{MR2504947} and \cite{MR2811131}, a characterization of those cubelike graphs  that admit perfect state transfer at a given time (for example $t = \frac{\pi}{8}$) is still unknown. 

In \cite{chan2013}, Chan gives examples of graphs with perfect state transfer which are constructed from unions of classes in the binary Hamming scheme. The determinations of which unions of distance classes in metric association schemes admit perfect state transfer would be an interesting problem. 

\bibliographystyle{plain}

\end{document}